\documentclass[letterpaper,10pt]{article} 

\usepackage{opticameet3} 
\usepackage{amsmath,amssymb} 
\usepackage[colorlinks=true,bookmarks=false,citecolor=blue,urlcolor=blue]{hyperref} 
\graphicspath{{Images/}}
\begin{document}

\title{New presences of $\pi$ and $e$ in Pascal's triangle}

\author{Mauricio Guevara Valerio}
\address{University of Costa Rica}
\email{mauricio.guevaravalerio@ucr.ac.cr}

\copyrightyear{2023}

\begin{abstract}
The following work shows new connections between the constants $\pi$ and $e$ with Pascal's triangle and the Lucas triangle, established via Fibonacci polynomials and similar means. Furthermore,  relations between the two famous constants and the rows of Pascal's triangle and the Lucas triangle are conjectured, together with some other important related identities. 
\end{abstract} 

\section{Background}
Discoveries of  $e$ and $\pi$ in Pascal's triangle are relatively recent in the history of mathematics. In the case of $e$, Harlan J. Brothers found the following remarkable relation in 2012 [1]:\\\medskip
\medskip
$$\lim_{n \to \infty} \dfrac{\frac{s_{n+1}}{s_n}} {\frac{s_n}{s_{n-1}}} = e$$
\\\medskip
\medskip
Where $s_n$ is the product of all the entries in the $n^{th}$ row of Pascal's triangle. 
\\\medskip
It is less known, however, that the great popularizer of mathematics Martin Gardner had already mentioned a relation between $e$ and the triangle through the Fibonacci numbers [2], namely the following (although as he stated it the claim is false):\medskip\\
\begin{center}
$e = \dfrac{1 + 1 + \frac{2}{2!} + \frac{3}{3!} + \frac{5}{4!} + \frac{8}{5!} + \frac{13}{6!} + \frac{21}{7!} + \frac{34}{8!} + \frac{55}{9!} + \cdots}{1 + 0 + \frac{1}{2!} + \frac{1}{3!} + \frac{2}{4!} + \frac{3}{5!} + \frac{5}{6!} + \frac{8}{7!} + \frac{13}{8!} + \frac{21}{9!} + \cdots}
$
\medskip\\
\end{center}
A slight alteration of the signs is all it takes to get the true identity:\medskip\\
\begin{center}
$e = \dfrac{1 + 1 + \frac{2}{2!} + \frac{3}{3!} + \frac{5}{4!} + \frac{8}{5!} + \frac{13}{6!} + \frac{21}{7!} + \frac{34}{8!} + \frac{55}{9!} + \cdots}{1 - 0 + \frac{1}{2!} - \frac{1}{3!} + \frac{2}{4!} - \frac{3}{5!} + \frac{5}{6!} - \frac{8}{7!} + \frac{13}{8!} - \frac{21}{9!} + \cdots}$ 
\end{center}
Gardner does not mention who discovered this curious identity, nor how to prove it. He — or rather his character “O'Shea” — says only that it was found “on the Web”. \\
Although this identity is very remarkable on its own, since it shows the connection between $e$ and the Fibonacci numbers, it turns out there is a similar and more interesting relationship between $e$ and the Fibonacci polynomials, through which a relationship between $e$, $\pi$ and Pascal's triangle can be established . \medskip
The structure of this work is as follows: in section 2 it will be shown how the identity mentioned by Gardner can be proved after the signs are fixed, using a more general formula relating $e$ and Fibonacci numbers, as well as other Lucas sequences. In section 3 a similar relation between $e$ and the Fibonacci Polynomials is established. In section 4 the connection with $\pi$ is established. Section 5 consists of several conjectures related to identities involving the two constant, Pascal's triangle and the Lucas triangle. Finally, in section 6 some other important related identities are demonstrated, together with 2 other conjectures related to the Lucas triangle. \medskip\\
\section{Identities relating e to Fibonacci numbers and other Lucas sequences}
\medskip
\flushleft\textbf{Theorem 1:}\\\medskip
\centering
$e=\dfrac{\sum\limits_{k=0}^\infty\frac{F_{k+1}+ xF_{k-1}}{k!}}{\sum\limits_{k=0}^\infty{\left(-1\right)^{k}\cdot \frac{F_{k-1}+xF_{k+1}}{k!}}}$\\\medskip 
\flushleft \textit{Proof:} From the Euler - Binet formula:\\\medskip
\centering$F_n=\dfrac{\phi^n-(1 -\phi) ^{n}}{\sqrt{5}}$\\\medskip 
\flushleft And the definition of $e^x$ as the sum of an infinite series:\\\medskip
\centering $e^x=\sum\limits_{k=0}^\infty \frac{x^k}{k!}$\\\medskip
\flushleft It follows that:\\\medskip
$$\dfrac{\sum\limits_{k=0}^\infty\frac{F_{k+1}+ xF_{k-1}}{k!}}{\sum\limits_{k=0}^\infty{\left(-1\right)^{k}\cdot \frac{F_{k-1}+xF_{k+1}}{k!}}}$$
\\\medskip
$$=\dfrac{\sum\limits_{k=0}^\infty\frac{\phi^{k+1}-(1-\phi)^{k+1}+x\phi^{k-1}-x(1-\phi)^{k-1}}{\sqrt{5}\cdot k!}}{\sum\limits_{k=0}^\infty(-1) ^k\cdot\frac{\phi^{k-1}-(1-\phi)^{k-1}+x\phi^{k+1}-x(1-\phi)^{k+1}}{\sqrt{5}\cdot k!}}$$
\\\medskip
It can be verified [3] that this last expression is equal to:\\\medskip
$$=\dfrac{\frac{2e^{1-\phi}(x\sqrt{5}e^{\sqrt{5}}-xe^{\sqrt{5}}+x\sqrt{5}+x+\sqrt{5}e^{\sqrt{5}}+e^{\sqrt{5}}+\sqrt{5}-1)}{4\sqrt{5}}}{\frac{2e^{-\phi}(x\sqrt{5}e^{\sqrt{5}}-xe^{\sqrt{5}}+x\sqrt{5}+x+\sqrt{5}e^{\sqrt{5}}+e^{\sqrt{5}}+\sqrt{5}-1)}{4\sqrt{5}}}$$\\\medskip
From here the identity follows almost immediately:\\\medskip
$$=\dfrac{2e^{1-\phi}}{2e^{-\phi}}=\dfrac{e^{1-\phi}}{e^{-\phi}}=e^{1-\phi+\phi}=e$$\\\medskip
And so we have proved our first Theorem. To get the identity Gardner mentions, we just have to set $x=0$ in Theorem 1:\\\medskip
$$e=\dfrac{\sum\limits_{k=0}^\infty\frac{F_{k+1}} {k!}}{\sum\limits_{k=0}^\infty{\left(-1\right)^{k}\cdot \frac{F_{k-1}} {k!}}}$$
If we set x = 1, we get a relation between $e$ and the Lucas numbers:\\\medskip
$$e = \dfrac{\frac{2}{0!}+ \frac{1}{1!} + \frac{3}{2!} + \frac{4}{3!} + \frac{7}{4!} + \frac{11}{5!} + \frac{18}{6!} + \frac{29}{7!} + \frac{47}{8!} + \frac{76}{9!} + \cdots}{\frac{2}{0!} - \frac{1}{1!} + \frac{3}{2!} - \frac{4}{3!} + \frac{7}{4!} - \frac{11}{5!} + \frac{18}{6!} - \frac{29}{7!} + \frac{47}{8!} - \frac{76}{9!} + \cdots}$$\\\medskip 
This identity can also be proved through a different formula, as may be seen in the final section. It is also worth mentioning that Theorem 1 gives relations between e and Lucas sequences that also have the same recurrence relation and link to the golden ratio, when x has positive integer values. For instance, setting x=2 gives the following identity:\\\medskip
$$e = \dfrac{\frac{3}{0!}+ \frac{1}{1!} + \frac{4}{2!} + \frac{5}{3!} + \frac{9}{4!} + \frac{14}{5!} + \frac{23}{6!} + \frac{37}{7!} + \frac{60}{8!} + \frac{97}{9!} + \cdots}{\frac{3}{0!} - \frac{2}{1!} + \frac{5}{2!} - \frac{7}{3!} + \frac{12}{4!} - \frac{19}{5!} + \frac{31}{6!} - \frac{50}{7!} + \frac{81}{8!} - \frac{131}{9!} + \cdots}$$\\\medskip
It may also be worth noting that by writing down the terms with negative indices of the sequence $3, 1,4, 5, 9,\cdots$ to the left, we can observe the following pattern:\\\medskip
$\cdots, -19, 12,-7, 5,-2,3,1,4,5, 9,\cdots$\\\medskip 
This pattern seems to hold for the other Lucas sequences as well. \\\medskip 
\section{The connection between $e$ and Pascal's triangle via Fibonacci polynomials} To prove the identity that relates Euler's number with Pascal's triangle through Fibonacci polynomials, we will use the closed form formula for the Fibonacci polynomials [4] instead of Binet's formula:\medskip\\
\begin{center}
$F_k(x)=\dfrac{(x+\sqrt{x^2+4})^k-(x-\sqrt{x^2+4})^k}{2^{k}\sqrt{x^2+4}}$  
\end{center}
\textbf{Theorem 2:}\medskip
\begin{center}
$e^x=\dfrac{\sum\limits_{k=1}^\infty\frac{F_k(x)}{k!}}{\sum\limits_{k=1}^\infty{\left(-1\right)^{k+1}\cdot \frac{F_k(x)}{k!}}}
$
\end{center}
\textit{Proof:}\hspace{3mm}From the closed form formula for the Fibonacci polynomials, it follows that:\\\medskip
\begin{center}
$\dfrac{\sum\limits_{k=1}^\infty\frac{F_k(x)}{k!}}{\sum\limits_{k=1}^\infty{\left(-1\right)^{k+1}\cdot \frac{F_k(x)}{k!}}}=\dfrac{\sum\limits_{k=1}^\infty\frac{\frac{\left(x+\sqrt{x^2+4}\right)^k -\left(x-\sqrt{x^2+4}\right)^k}{2^k\sqrt{x^2+4}}}{k!}}{\sum\limits_{k=1}^\infty\left(-1\right)^{k+1}\frac{\frac{\left(x+\sqrt{x^2+4}\right)^k-\left(x-\sqrt{x^2+4}\right)^k}{2^k\sqrt{x^2+4}}}{k!}}$
\end{center}
\medskip
\medskip
\begin{center}
$=\dfrac{\frac{1}{\sqrt{x^2+4}}\cdot\sum\limits_{k=1}^\infty\dfrac{\left(\frac{x+\sqrt{x^2+4}}{2}\right)^k-\left(\frac{x-\sqrt{x^2+4}}{2}\right)^k}{k!}}{\frac{1}{\sqrt{x^2+4}}\cdot\sum\limits_{k=1}^\infty\left(-1\right)^{k+1}\dfrac{\left(\frac{x+\sqrt{x^2+4}}{2}\right)^k-\left(\frac{x-\sqrt{x^2+4}}{2}\right)^k}{k!}}$
\end{center}
\medskip
From the definition of $e^x$ as an infinite series, it follows that:\\\medskip
\begin{center}
$=\dfrac{e^{\frac{x+\sqrt{x^2+4}}{2}}-e^{\frac{x-\sqrt{x^2+4}}{2}}}{\sum\limits_{k=1}^\infty\dfrac{-\left(-\frac{x+\sqrt{x^2+4}}{2}\right)^k--\left(-\frac{x-\sqrt{x^2+4}}{2}\right)^k}{k!}}$
\medskip
\medskip
$=\dfrac{e^{\frac{x+\sqrt{x^2+4}}{2}}-e^{\frac{x-\sqrt{x^2+4}}{2}}}{\sum\limits_{k=1}^\infty\dfrac{-\left(-\frac{x+\sqrt{x^2+4}}{2}\right)^k+\left(-\frac{x-\sqrt{x^2+4}}{2}\right)^k}{k!}}$
\medskip
$=\dfrac{e^{\frac{x+\sqrt{x^2+4}}{2}}-e^{\frac{x-\sqrt{x^2+4}}{2}}} {e^{-\frac{x-\sqrt{x^2+4}}{2}}-e^{-\frac{x+\sqrt{x^2+4}}{2}}}$
\medskip
$=\dfrac{e^{\frac{x+\sqrt{x^2+4}}{2}}-e^{\frac{x-\sqrt{x^2+4}}{2}}} {e^{\frac{-x+\sqrt{x^2+4}}{2}}-e^{\frac{-x-\sqrt{x^2+4}}{2}}}$
\medskip
\end{center}
\medskip
Now we just have to simplify a bit:
\\\medskip
\begin{center}
$=\dfrac{\left(e^{\sqrt{x^2+4}}-1\right)\cdot e^{\frac{x-\sqrt{x^2+4}}{2}}}{\left(e^{\sqrt{x^2+4}}-1\right)\cdot e^{\frac{-x-\sqrt{x^2+4}}{2}}}$\medskip
$=\dfrac{e^{\frac{x-\sqrt{x^2+4}}{2}}}{e^{\frac{-x-\sqrt{x^2+4}}{2}}}$\medskip
$=e^{\frac{x-\sqrt{x^2+4}}{2}}\cdot e^{\frac{x+\sqrt{x^2+4}}{2}}$
\medskip
$=e^{\frac{2x}{2}}$
\medskip
$=e^x$
\medskip
\end{center}
Thus completing the proof. To get to an identity very similar to the one mentioned by Gardner, we just have to set $x=1$ in Theorem 2: \\\medskip
\begin{center}
$e = \dfrac{1 + \frac{1}{2!} + \frac{2}{3!} + \frac{3}{4!} + \frac{5}{5!} + \frac{8}{6!} + \frac{13}{7!} + \frac{21}{8!} + \frac{34}{9!} + \cdots}{1 - \frac{1}{2!} + \frac{2}{3!} - \frac{3}{4!} + \frac{5}{5!} - \frac{8}{6!} + \frac{13}{7!} - \frac{21}{8!} + \frac{34}{9!} - \cdots}$
\medskip
\end{center}
And it is known that the Fibonacci polynomials can be found in the shallow diagonals of Pascal's triangle, if we interpret its numbers as coefficients:\\\medskip
\begin{figure}[htp]
    \centering
   \includegraphics[width=4cm]{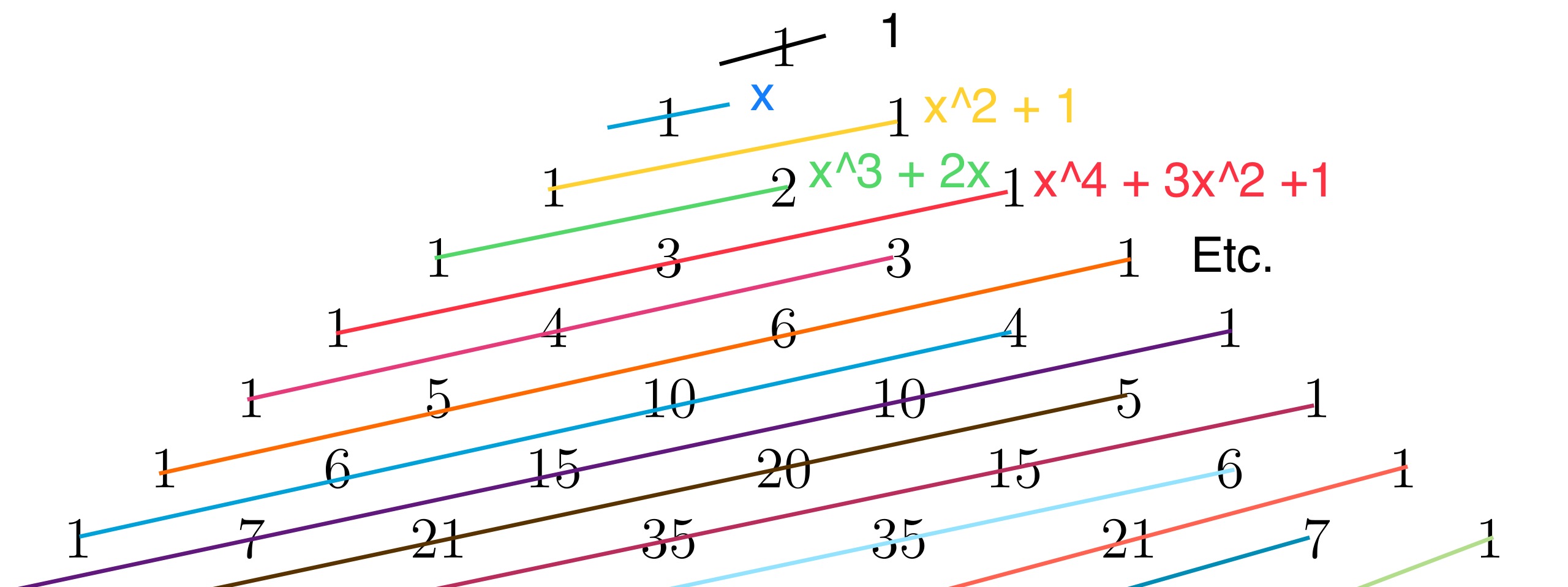}
    \caption{Fibonacci polynomials in Pascal's triangle}
    \label{fig:PascalTriangleFibonacci}
\end{figure}
\section{The connection between $\pi,e$ and Pascal's triangle via Fibonacci polynomials}
\textbf{Theorem 3:}
\medskip
\begin{center}
$$\frac{e-\frac{1}{e}}{2}=\frac{\pi^2}{3!}-\frac{\pi^4-3\pi^2}{5!}+\frac{\pi^6-5\pi^4+6\pi^2}{7!}-\frac{\pi^8-7\pi^6+15\pi^4-10\pi^2}{9!}+\cdots$$
\end{center}
\medskip
\textit{Proof:}\hspace{2mm}From Euler's identity:\medskip
\begin{center}
$e^{i\pi}+1=0\hspace{6mm}\rightarrow \hspace{6mm} e^{i\pi}=-1$
\end{center}
Setting $x=i\pi$ in Theorem 2 yields:\medskip
\begin{center}
$e^{i\pi} = \dfrac{1 + \frac{i\pi}{2!} + \frac{(i\pi)^2+1}{3!} + \frac{(i\pi)^3+2i\pi}{4!} + \frac{(i\pi)^4+3(i\pi)^2+1}{5!} + \frac{(i\pi)^5+4(i\pi)^3+3i\pi}{6!} + \cdots}{1 - \frac{i\pi}{2!} + \frac{(i\pi)^2+1}{3!} - \frac{(i\pi)^3+2i\pi}{4!} + \frac{(i\pi)^4+3(i\pi)^2+1}{5!} - \frac{(i\pi)^5+4(i\pi)^3+3i\pi}{6!} + \cdots}$
\end{center}
\medskip
Combining the above with Euler's identity gives the following identity:
\begin{center}
\medskip
\medskip
$-1=\dfrac{1 + \frac{i\pi}{2!} + \frac{(i\pi)^2+1}{3!} + \frac{(i\pi)^3+2i\pi}{4!} + \frac{(i\pi)^4+3(i\pi)^2+1}{5!} + \frac{(i\pi)^5+4(i\pi)^3+3i\pi}{6!} + \cdots}{1 - \frac{i\pi}{2!} + \frac{(i\pi)^2+1}{3!} - \frac{(i\pi)^3+2i\pi}{4!} + \frac{(i\pi)^4+3(i\pi)^2+1}{5!} - \frac{(i\pi)^5+4(i\pi)^3+3i\pi}{6!} + \cdots}$\medskip

$\rightarrow-1 + \frac{i\pi}{2!} - \frac{(i\pi)^2+1}{3!} + \frac{(i\pi)^3+2i\pi}{4!} - \frac{(i\pi)^4+3(i\pi)^2+1}{5!} + \frac{(i\pi)^5+4(i\pi)^3+3i\pi}{6!} - \cdots=1 + \frac{i\pi}{2!} + \frac{(i\pi)^2+1}{3!} + \frac{(i\pi)^3+2i\pi}{4!} + \frac{(i\pi)^4+3(i\pi)^2+1}{5!} + \frac{(i\pi)^5+4(i\pi)^3+3i\pi}{6!} + \cdots$
\end{center}
\medskip
Now all the terms whose denominator is an even factorial get cancelled out:\medskip
\begin{center}
$\rightarrow-1 - \frac{(i\pi)^2+1}{3!} - \frac{(i\pi)^4+3(i\pi)^2+1}{5!} - \cdots =1 + \frac{(i\pi)^2+1}{3!} + \frac{(i\pi)^4+3(i\pi)^2+1}{5!} + \cdots$
\end{center}
\medskip
Since we have an identity of the form $-x=x$, that means both sides of the equation must equal $0$. We can therefore take the right hand side:\medskip
\begin{center}
$0=1 + \frac{(i\pi)^2+1}{3!} + \frac{(i\pi)^4+3(i\pi)^2+1}{5!} + \cdots$\medskip
$=1 - \frac{\pi^2-1}{3!} + \frac{\pi^4-3\pi^2+1}{5!} -\frac{\pi^6-5\pi^4+6\pi^2-1}{7!}+ \frac{\pi^8-7\pi^6+15\pi^4-10\pi^2+1}{9!}-\cdots$
\end{center}
\medskip
And use the Taylor series expansion for $\sinh(1) $:\medskip
\begin{center}
$sinh(1)=\frac{e-\frac{1}{e}}{2}= 1+\frac{1}{3!}+\frac{1}{5!}+\frac{1}{7!}+\frac{1}{9!}+\cdots=1-(-\frac{1}{3!})+\frac{1}{5!}-(-\frac{1}{7!})+\frac{1}{9!}-\cdots$
\end{center}
\medskip
By combining the two previous identities and simplifying, we get:\medskip
\begin{center}
$0=\frac{e-\frac{1}{e}}{2}- \frac{\pi^2}{3!} + \frac{\pi^4-3\pi^2}{5!} -\frac{\pi^6-5\pi^4+6\pi^2}{7!}+ \frac{\pi^8-7\pi^6+15\pi^4-10\pi^2}{9!}-\cdots$
\end{center}
\medskip
And finally we get the beautiful identity that relates $e,\pi$ and Pascal's triangle through Fibonacci polynomials, completing the proof of Theorem 3:\medskip
\begin{center}
$\frac{e-\frac{1}{e}}{2}=\frac{\pi^2}{3!}-\frac{\pi^4-3\pi^2}{5!}+\frac{\pi^6-5\pi^4+6\pi^2}{7!}-\frac{\pi^8-7\pi^6+15\pi^4-10\pi^2}{9!}+\cdots$\medskip
\end{center}
The pattern in the coefficients of the powers of $\pi$ can be found in the ``rotated'' Pascal triangle, following the lines in the second image\medskip
\begin{figure}[htp]
    \centering
    \includegraphics[width=4cm]{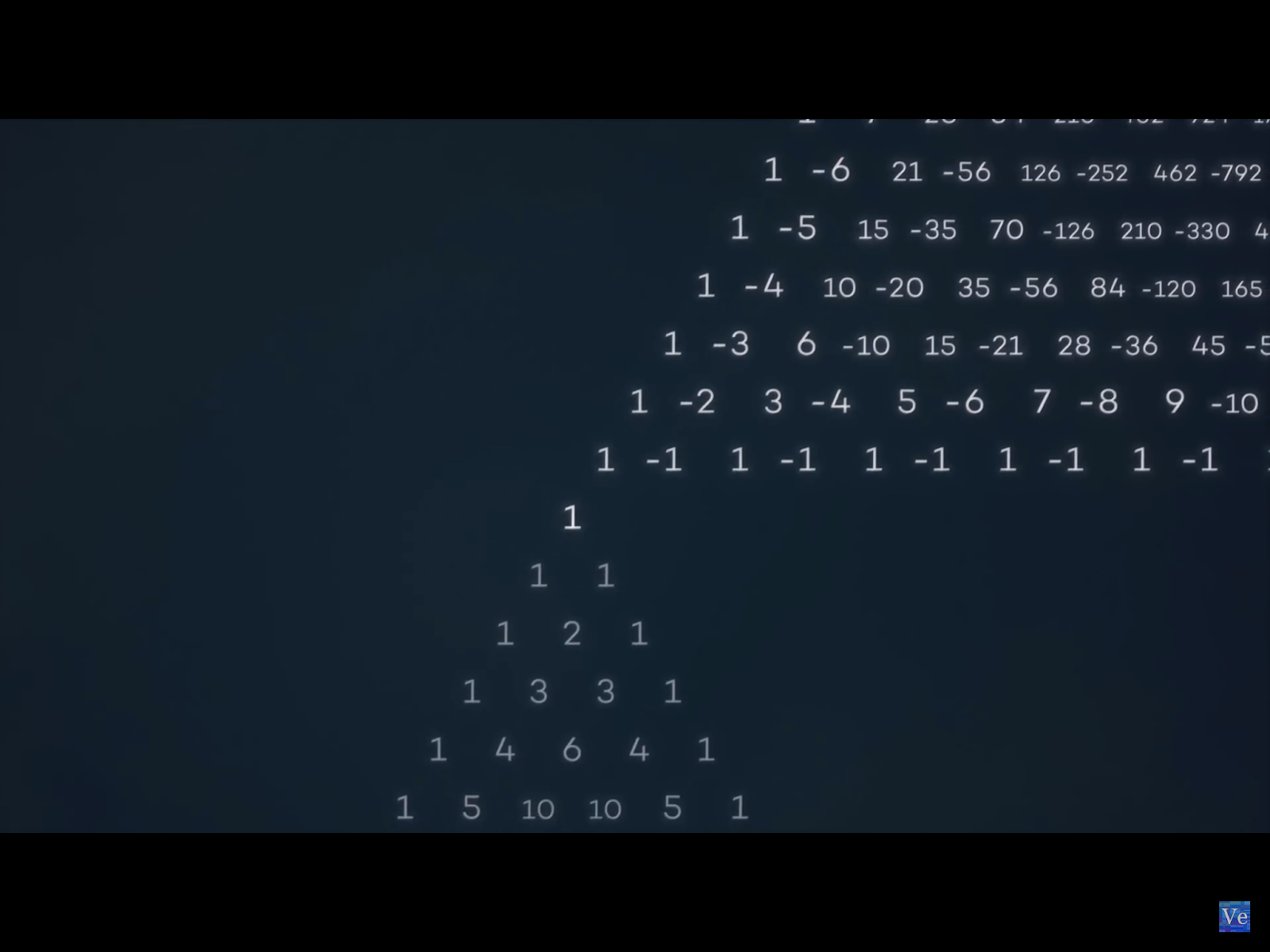}
    \caption{Extension of Pascal's triangle to negative numbers: the “rotated” Pascal triangle}
    \label{fig:Pascal}
\end{figure}
\\
\begin{figure}[htp]
    \centering
   \scalebox{0.35}{\includegraphics[width=4cm]{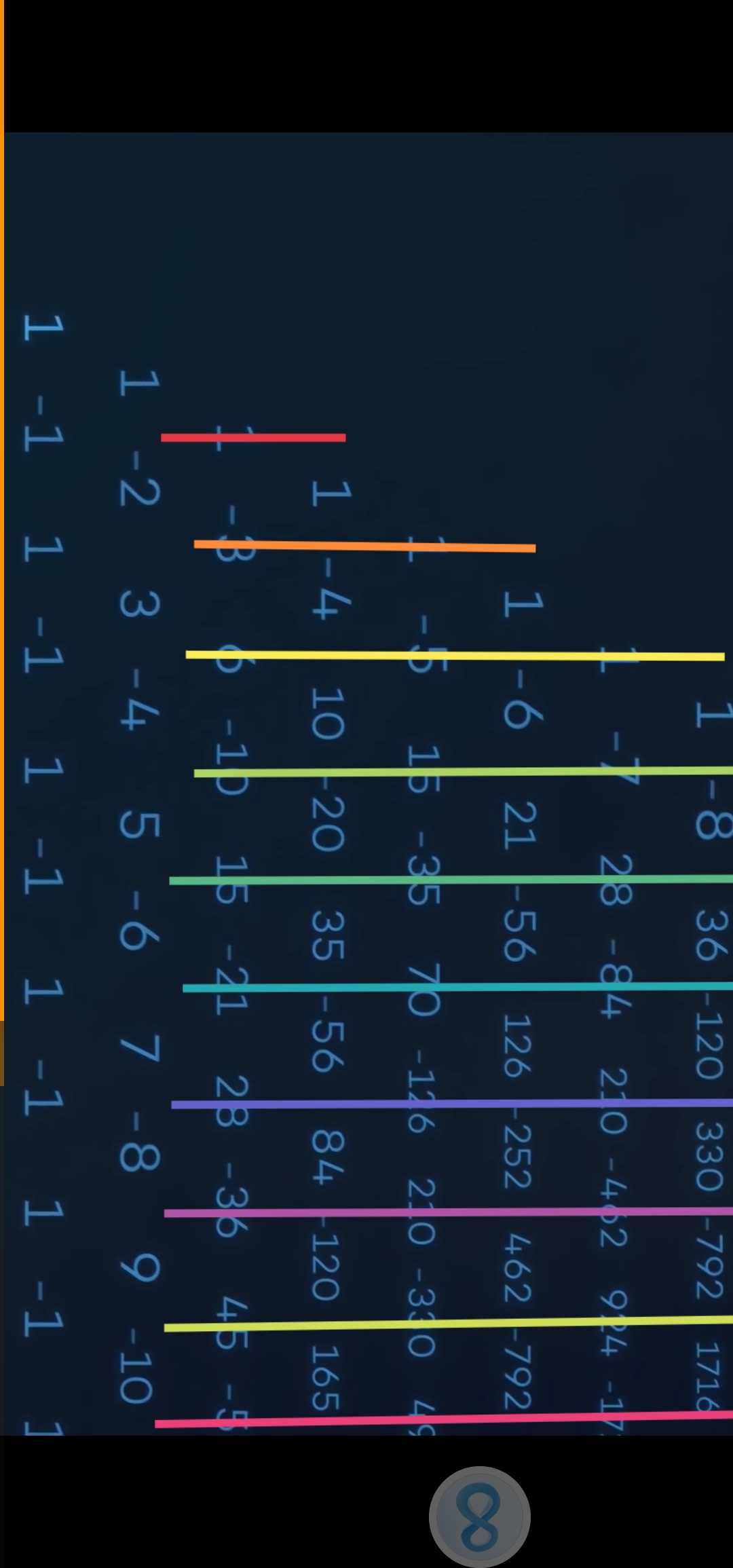}}
    \caption{Pattern of the coefficients for the powers of $\pi$ in Theorem $3$}
    \label{fig:Rotated}
\end{figure}
\section{Conjectures}
A possible way of proving Theorem 3 without using Euler's identity is the following:\\
\medskip
\medskip
\medskip
\textbf{Conjecture 1.} \\\medskip 
$\frac{1}{1!}=\frac{\pi^2}{3!}-\frac{\pi^4}{5!}+\frac{\pi^6}{7!}-\frac{\pi^8}{9!}+\cdots$\\\medskip

$\frac{1}{3!}=\frac{3\pi^2}{5!}-\frac{5\pi^4}{7!}+\frac{7\pi^6}{9!}-\frac{9\pi^8}{11!}+\cdots$\\\medskip

$\frac{1}{5!}=\frac{6\pi^2}{7!}-\frac{15\pi^4}{9!}+\frac{28\pi^6}{11!}-\frac{45\pi^8}{13!}+\cdots$\\\medskip

$\frac{1}{7!}=\frac{10\pi^2}{9!}-\frac{35\pi^4}{11!}+\frac{84\pi^6}{13!}-\frac{165\pi^8}{15!}+\cdots$\\\medskip
And so on ad infinitum. Using Wolfram Alpha to verify these identities, it can be seen that this pattern seems to hold. If this is true, then we can add the terms nicely like this: \\\medskip
\begin{figure}
\centering
   \scalebox{0.2}{\includegraphics{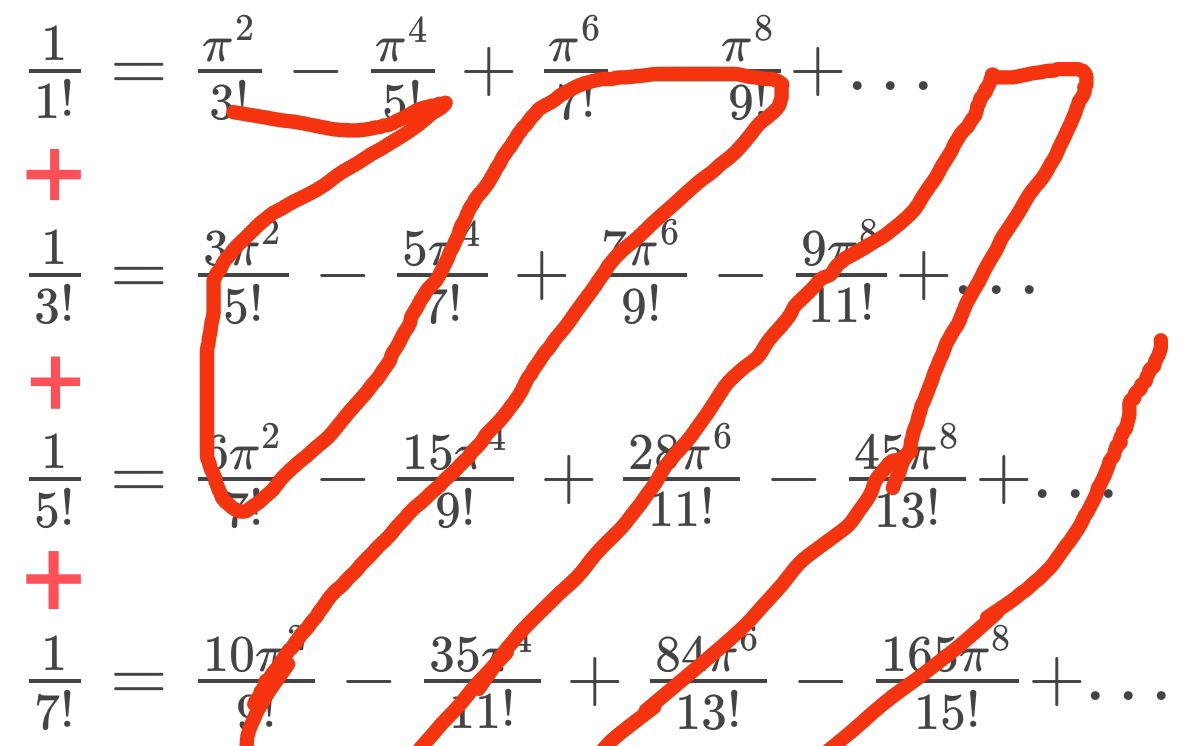}}
    \caption{Zigzag proof of Theorem 3}
    \label{fig:my_label}
\end{figure}
And rather surprisingly, there seems to be another beautiful relation between the 2 constants and Pascal's triangle in a single identity. The identity is the following:\\\medskip
\textbf{Conjecture 2.}
$$e=\frac{\pi^2}{3!}+\frac{4\pi^2}{4!}-\frac{\pi^4-10\pi^2}{5!}-\frac{6\pi^4-20\pi^2}{6!}+\frac{\pi^6-21\pi^4+35\pi^2}{7!}+\frac{8\pi^6-56\pi^4+56\pi^2}{8!}-\frac{\pi^8-36\pi^6+126\pi^4-84\pi^2}{9!}-\cdots$$\\\medskip
\begin{figure}[ht]
    \centering
  \scalebox{0.3} {\includegraphics{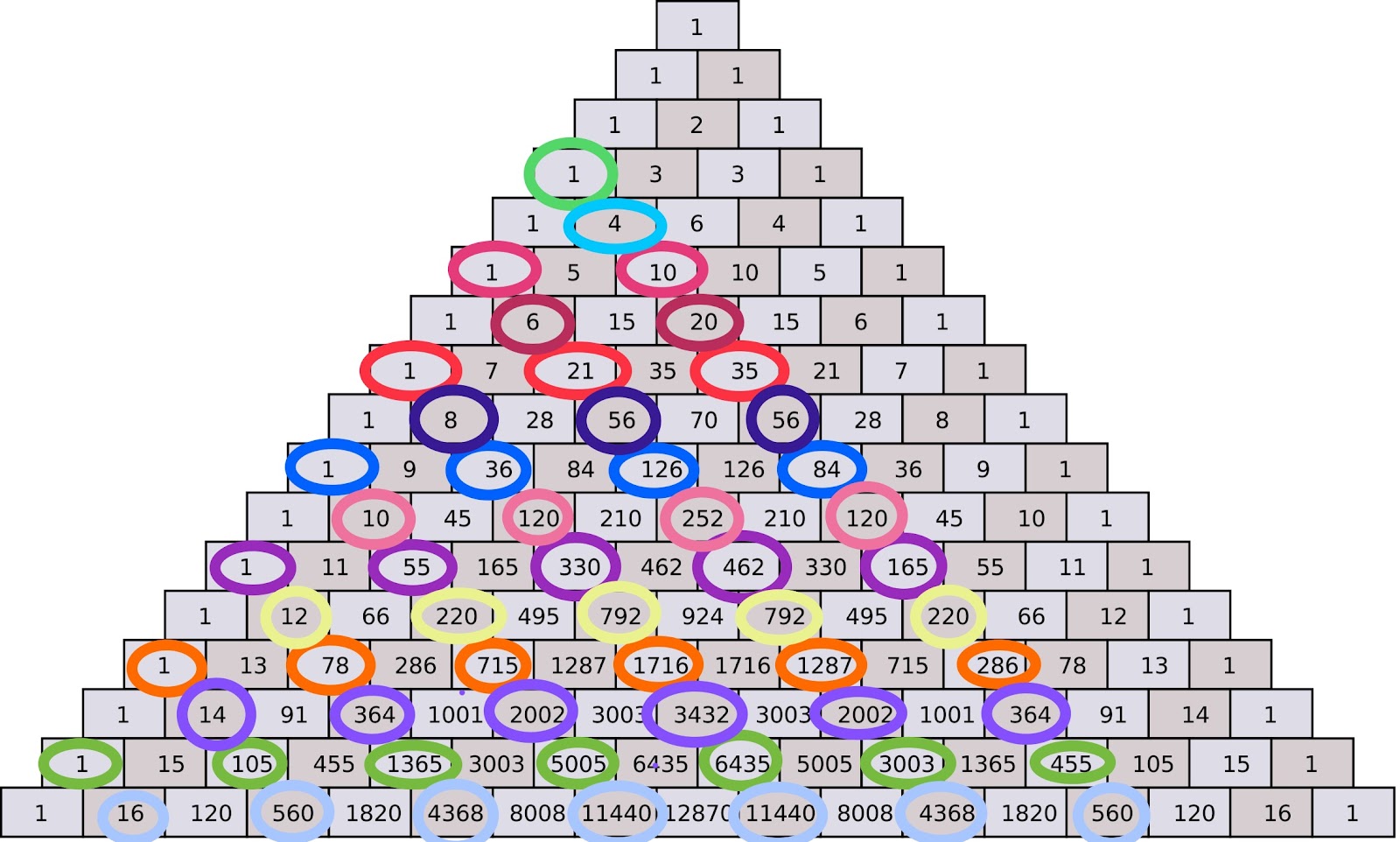} }
    \caption{Pattern in Conjecture 2}
    \label{fig:my_label}
\end{figure}
A good reason for believing this is true, besides empirical computation, is a similar “zig zag” representation for the following infinite series expression for $e$:\\\medskip
$$e=\frac{2}{1!}+\frac{4}{3!}+\frac{6}{5!}+\frac{8}{7!}+\frac{10}{9!}+\cdots$$\\\medskip 
\begin{figure}[ht]
    \centering
    \scalebox{0.2}{\includegraphics{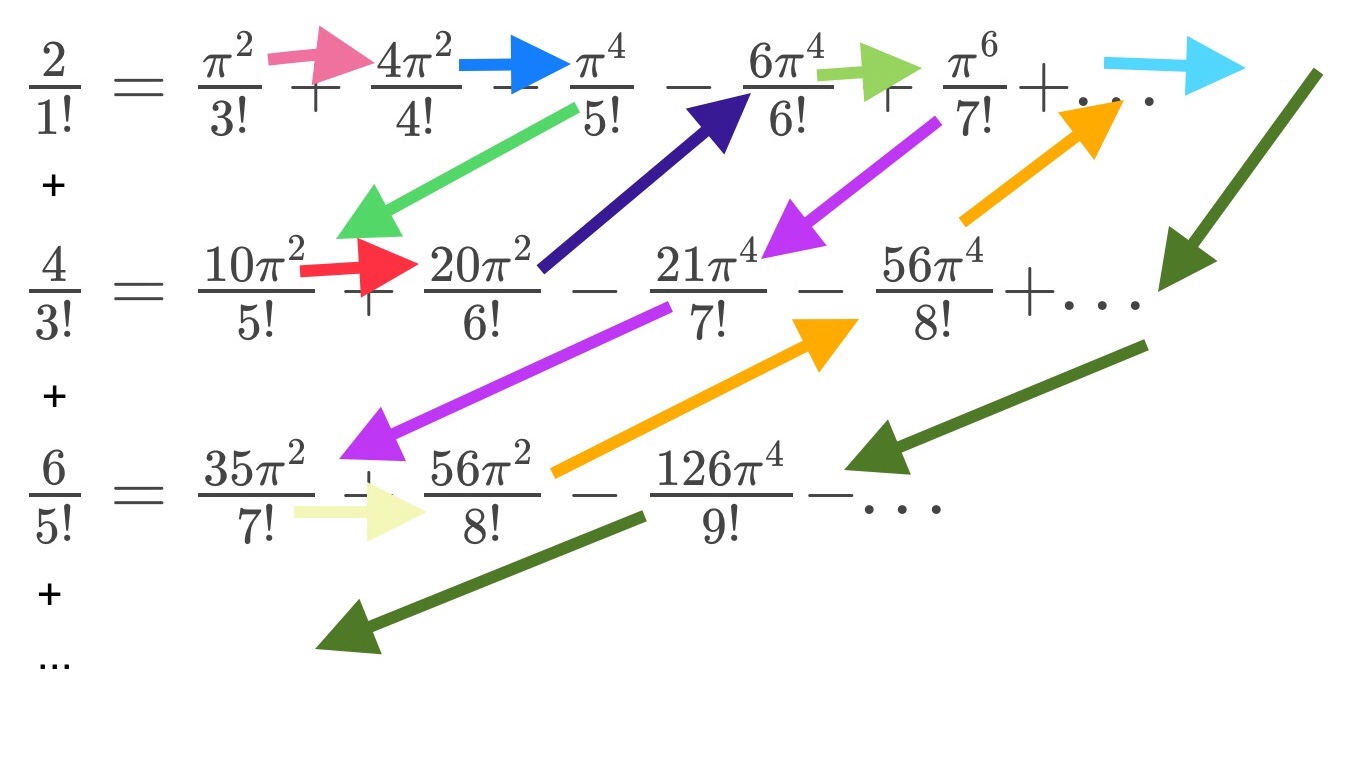}}
    \caption{Zig zag proof of Conjecture 2}
    \label{fig:my_label}
\end{figure}
\newpage
Another reason is that one can find similar evidence that if we split the sum into two infinite series, they seem to add up respectively to the hyperbolic sine of 1 and the hyperbolic cosine of 1, and adding them both would give us the infinite series for e if they had those values: \\\medskip
\textbf{Conjecture 3.} $$\frac{e-\frac{1}{e}}{2}=\frac{4\pi^2}{4!}-\frac{6\pi^4-20\pi^2}{6!}+\frac{8\pi^6-56\pi^4+56\pi^2}{8!}-\frac{10\pi^8-120\pi^6+252\pi^4-120\pi^2}{10!}+\cdots$$
\\\medskip
\textbf{Conjecture 4.}$$\frac{e+\frac{1}{e}}{2}=\frac{\pi^2}{3!}-\frac{\pi^4-10\pi^2}{5!}+\frac{\pi^6-21\pi^4+35\pi^2}{7!}-\frac{\pi^8-36\pi^6+126\pi^4-84\pi^2}{9!}+\cdots$$\\\medskip

Because of a notable property of Pascal's triangle, Conjecture 2 also shows an interesting relation between $e$, $\pi$ and the prime numbers: in all the terms of the sum whose denominator is the factorial of a prime number, all the terms in the numerator (with the exception of the first term in each numerator) are multiples of that prime number, multiplied by even powers of $\pi$. 
\\\medskip 
We can use this fact to get a series for e involving only the factorials of composite numbers in the denominator: since $2 = \frac{\pi^2}{3!} + \frac{4\pi^2}{4!} - \frac{\pi^4}{5!} - \frac{6\pi^4}{6!}+\cdots$, we can get rid of all the terms with coefficient $1$ times some even power of $\pi$ in the numerator, by substituting that series with $2$, after which we can take the common prime factor and simplify the resulting expression using the obvious facts that $(pk)/p! = k/(p-1)!$ and that $p-1$ is never a prime number when $p$ is a prime number greater than 3, to get this identity:\\\medskip $$e=2+\frac{2\pi^2}{4!}-\frac{3\pi^4-25\pi^2}{6!}-\frac{56\pi^4-56\pi^2}{8!}+\frac{36\pi^6-126\pi^4+84\pi^2}{9!}-\frac{5\pi^8-150\pi^6+294\pi^4-135\pi^2}{10!}+\cdots$$
\\\medskip
\section{Other related identities}
In this section I will mention some other important identities that follow from or are similar to ones that have been already mentioned. \\\medskip
\medskip
A connection between e and $\pi$ established via Chebyshev polynomials of the second kind can be deduced setting $x=2\pi i$ in Theorem 2:\\
\medskip
$$\frac{e-\frac{1}{e}}{2}=\frac{8\pi^2}{4!}-\frac{32\pi^4-32\pi^2}{6!}+\frac{128\pi^6-192\pi^4+80\pi^2}{8!}-\frac{512\pi^8-1024\pi^6+672\pi^4-160\pi^2}{10!}+ \cdots$$\\\medskip

The proof of this identity is very analogous to the proof of Theorem 3. The only significant difference is that since $e^{2\pi i}=1$, the terms that get canceled out are the ones whose denominator is the factorial of an odd number. \\\medskip 

It is also worth noting that many other interesting identities can be obtained through a similar connection between $e$ and the Lucas polynomials, through which a connection between $e, \pi$ and the Lucas triangle can be established:\\\medskip 

\textbf{Theorem 4. 
} 
\begin{center}
$e^x=\dfrac{\sum\limits_{k=0}^\infty\frac{L_k(x)}{k!}}{\sum\limits_{k=0}^\infty{\left(-1\right)^{k}\cdot \frac{L_k(x)}{k!}}}
$ 
\end{center}
\textit{Proof:}\hspace{3mm} From the closed form formula for the Lucas polynomials [5]:\\\medskip

$$L_k(x)=\Big(\frac{1}{2}\Big)^{k}\cdot \Big[\Big(x+\sqrt{x^2+4}\Big)^k+\Big(x-\sqrt{x^2+4}\Big)^k\Big]$$\\\medskip
\medskip
And the definition of $e^x$ as an infinite series:\\\medskip

\centering$e^x=\sum\limits_{k=0}^\infty \frac{x^k}{k!}$ \\\medskip
\flushleft It follows that:\\\medskip 
$\dfrac{\sum\limits_{k=0}^\infty\frac{L_k(x)}{k!}}{\sum\limits_{k=0}^\infty{\left(-1\right)^{k}\cdot \frac{L_k(x)}{k!}}}=\dfrac{\sum\limits_{k=0}^\infty{\Big(\frac{1}{2}\Big)^k\frac{\Big[\Big(x+\sqrt{x^2+4}\Big)^k+\Big(x-\sqrt{x^2+4}\Big)^k\Big]}{k!}}} {\sum\limits_{k=0}^\infty{\frac{(-1)^k \cdot\Big(\frac{1}{2}\Big)^k \Big[\Big(x+\sqrt{x^2+4}\Big)^k+\Big(x-\sqrt{x^2+4}\Big)^k\Big]}{k!}}}$\\\medskip
$=\dfrac{\sum\limits_{k=0}^\infty{\frac{\Big(\frac{x+\sqrt{x^2+4}} {2}\Big)^k+\Big(\frac{x-\sqrt{x^2+4}} {2}\Big)^k}{k!}}}{\sum\limits_{k=0}^\infty{\frac{\Big(\frac{-\big(x+\sqrt{x^2+4}\big)} {2}\Big)^k+\Big(\frac{-\big(x-\sqrt{x^2+4}\big)}{2}\Big)^k}{k!}}}$
\\\medskip 
$=\dfrac{\sum\limits_{k=0}^\infty{\frac{\Big(\frac{x+\sqrt{x^2+4}} {2}\Big)^k+\Big(\frac{x-\sqrt{x^2+4}} {2}\Big)^k}{k!}}}{\sum\limits_{k=0}^\infty{\frac{\Big(\frac{-x-\sqrt{x^2+4}} {2}\Big)^k+\Big(\frac{-x+\sqrt{x^2+4}}{2}\Big)^k}{k!}}}$
\\\medskip
$=\dfrac{e^{\frac{x+\sqrt{x^2+4}}{2}}+e^{\frac{x-\sqrt{x^2+4}} {2}}} {e^{\frac{-x-\sqrt{x^2+4}}{2} }+e^{\frac{-x+\sqrt{x^2+4}}{2}}}$\\\medskip
\medskip
$=\dfrac{e^{\frac{x-\sqrt{x^2+4}}{2}} \Big(e^{\frac{2\sqrt{x^2+4}}{2}}+1 \Big) } {e^{\frac{-x-\sqrt{x^2+4}}{2}}\Big(e^{\frac{2\sqrt{x^2+4}}{2}}+1 \Big) } $
\\\medskip
$=\dfrac{e^{\frac{x-\sqrt{x^2+4}}{2}}}{e^{\frac{-x-\sqrt{x^2+4} }{2}}}$
\\\medskip
$=e^{\frac{x-\sqrt{x^2+4}}{2}} \cdot e^{\frac{x+\sqrt{x^2+4}}{2}} $
\\\medskip
$=e^{\frac{2x}{2}} $
\\\medskip
$=e^x$
\\\medskip 
As was to be proved. \\\medskip 
Setting $x=1$ gives again the identity for $e$ in terms of the Lucas numbers:\\\medskip 
$e = \dfrac{\frac{2}{0!}+ \frac{1}{1!} + \frac{3}{2!} + \frac{4}{3!} + \frac{7}{4!} + \frac{11}{5!} + \frac{18}{6!} + \frac{29}{7!} + \frac{47}{8!} + \frac{76}{9!} + \cdots}{\frac{2}{0!} - \frac{1}{1!} + \frac{3}{2!} - \frac{4}{3!} + \frac{7}{4!} - \frac{11}{5!} + \frac{18}{6!} - \frac{29}{7!} + \frac{47}{8!} - \frac{76}{9!} + \cdots}$\\\medskip 
\begin{figure}[ht]
    \centering
\scalebox{0.5}{\includegraphics{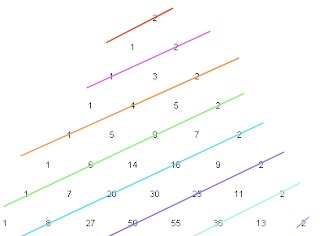}}
    \caption{Lucas polynomials in the shallow diagonals of the Lucas triangle}
    \label{fig:my_label}
\end{figure}
\newpage
Setting $x = i\pi$ gives the following identity after some straightforward simplification:\\\medskip 
$$e + e^{-1}=\frac{\pi^2}{2!} - \frac{\pi^4-4\pi^2}{4!}+ \frac{\pi^6 - 6\pi^4 + 9\pi^2}{6!} - \frac{\pi^8-8\pi^6+20\pi^4-16\pi^2}{8!} + \cdots$$\\\medskip
A more general relation between $e$ and the Lucas polynomials can be proved:\\\medskip 
\medskip
\textbf{Theorem 5:}
$$e^{L_n(x)}= \dfrac{{\sum\limits_{k=0}^\infty\frac{L_{nk} (x)}{k!}}} {{\sum\limits_{k=0}^\infty(-1)^k\frac{L_{nk} (x)}{k!}}}$$
\\\medskip 
\textit{Proof:} 
\\\medskip
$$\dfrac{{\sum\limits_{k=0}^\infty\frac{L_{nk} (x)}{k!}}}{{\sum\limits_{k=0}^\infty(-1)^k\frac{L_{nk} (x)}{k!}}}=\dfrac{\sum\limits_{k=0}^\infty\frac{\Big[\Big(\frac{x+\sqrt{x^2+4}}{2}\Big) ^n\Big]^k} {k!} +  \sum\limits_{k=0}^\infty\frac{\Big[\Big(\frac{x-\sqrt{x^2+4}}{2}\Big) ^n\Big]^k}{k! }}{\sum\limits_{k=0}^\infty\frac{\Big[-\Big(\frac{x+\sqrt{x^2+4}} {2}\Big) ^n\Big]^k} {k!}+\sum\limits_{k=0}^\infty\frac{\Big[-\Big(\frac{x-\sqrt{x^2+4}}{2}\Big)^n\Big]^k} {k!}} $$
\\\medskip
$$=\dfrac{e^{\Big(\frac{x+\sqrt{x^2+4}}{2}\Big) ^n}+ e^{\Big(\frac{x-\sqrt{x^2+4}} {2}\Big)^n}}  {\frac{1}{e^{\Big(\frac{x+\sqrt{x^2+4}}{2}\Big)^n}} + \frac{1}{e^{\Big(\frac{x-\sqrt{x^2+4}}{2}\Big)^n}}} $$
\\\medskip
$$=\dfrac{e^{\Big(\frac{x+\sqrt{x^2+4}}{2}\Big) ^n}+{e^{\Big(\frac{x-\sqrt{x^2+4}} {2}\Big)^n}}} {\frac{e^{\Big(\frac{x+\sqrt{x^2+4}}{2}\Big) ^n}+{e^{\Big(\frac{x-\sqrt{x^2+4}} {2}\Big)^n}}} {e^{{\Big(\frac{x+\sqrt{x^2+4}}{2}\Big) ^n}+\Big(\frac{x-\sqrt{x^2+4}} {2}\Big)^n}}}  $$
\\\medskip
$$=e^{{\Big(\frac{x+\sqrt{x^2+4}}{2}\Big) ^n}+\Big(\frac{x-\sqrt{x^2+4}} {2}\Big)^n}$$
\\\medskip 
$$=e^{L_n(x)}$$
\\\medskip 
As was to be proved. \\\medskip 
A couple of identities similar to the one previously conjectured based on the rows of Pascal's triangle seem to hold for the rows of the Lucas triangle as well:\\\medskip
\flushleft\textbf{Conjecture 5. }$$e=\frac{\pi^2}{2\cdot2!}-\frac{\pi^4-14\pi^2}{2\cdot{4!}}+\frac{\pi^6-27\pi^4+55\pi^2}{2\cdot6!}-\frac{\pi^8-44\pi^6+182\pi^4-140\pi^2}{2\cdot8!}+\cdots$$
\\\medskip
\flushleft\textbf{Conjecture 6.}\\\medskip
$$e=\frac{5\pi^2}{2\cdot3!}-\frac{7\pi^4-30\pi^2}{2\cdot5!}+\frac{9\pi^6-77\pi^4+91\pi^2}{2\cdot7!}-\frac{11\pi^8-156\pi^6+378\pi^4-204\pi^2}{2\cdot9!}+\cdots$$ 
\\\medskip 
Here once again we see a connection with the primes: in the numerator of all the terms whose denominator is 2 times the factorial of a prime number, all the terms except the first one are multiples of that prime number times even powers of $\pi$. 
\\\medskip
\begin{figure}[ht]
\scalebox{0.6}{\includegraphics{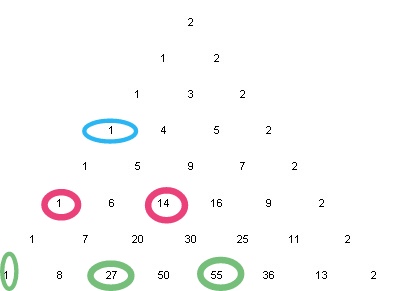}}
\centering
    {\caption{Pattern in Conjecture 5}
    \label{fig:my_label}}
\end{figure}
\begin{figure}[ht]
\centering
\scalebox{0.14}{\includegraphics{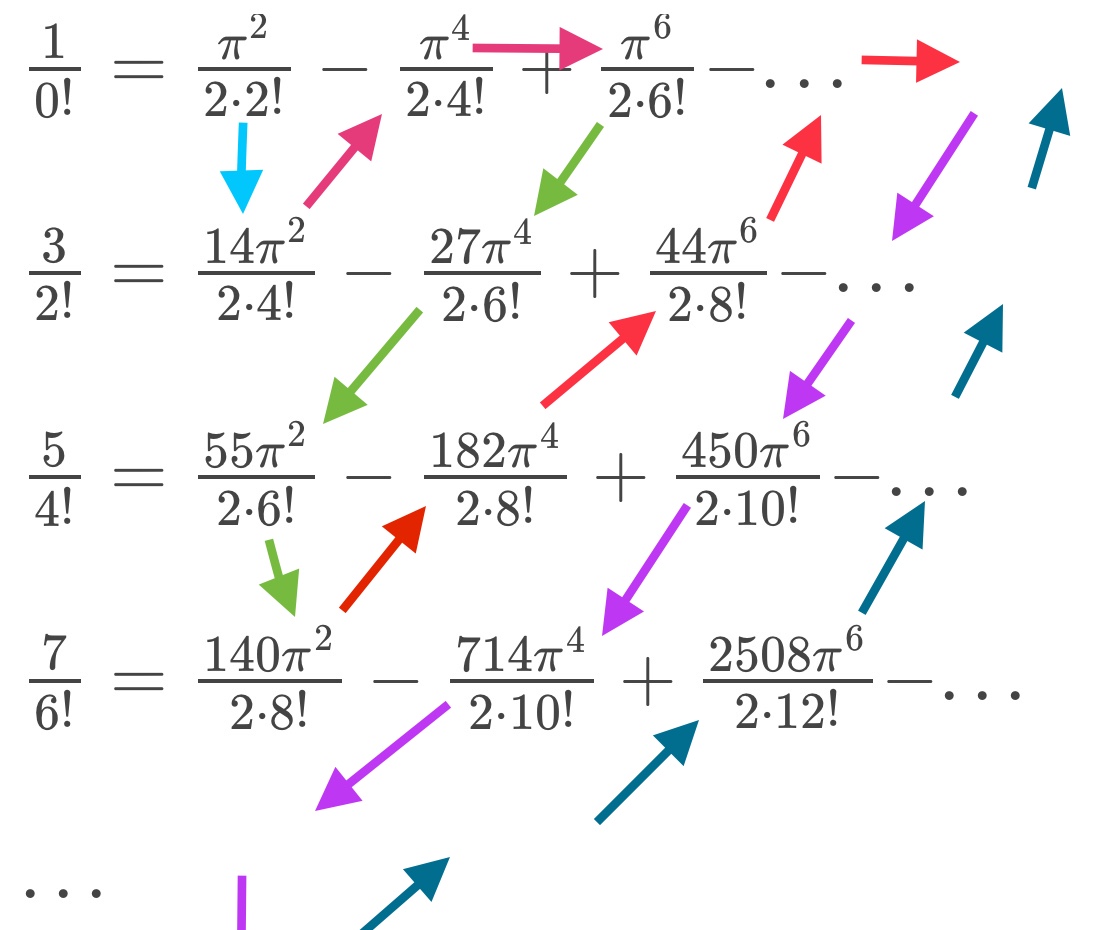}} 
\caption{Zig zag proof of Conjecture 5}
\label{fig:my_label}
\end{figure}
\newpage
\begin{figure}[ht]
    \centering
\scalebox{0.7}{\includegraphics{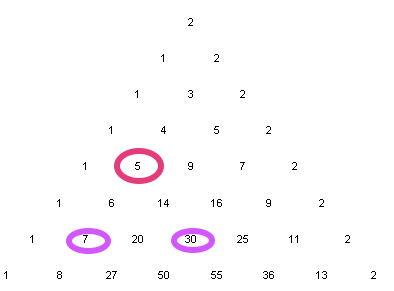}}
    \caption{Pattern in Conjecture 6}
    \label{fig:my_label}
\end{figure}
\begin{figure}[ht]
    \centering
    \scalebox{0.14}{\includegraphics{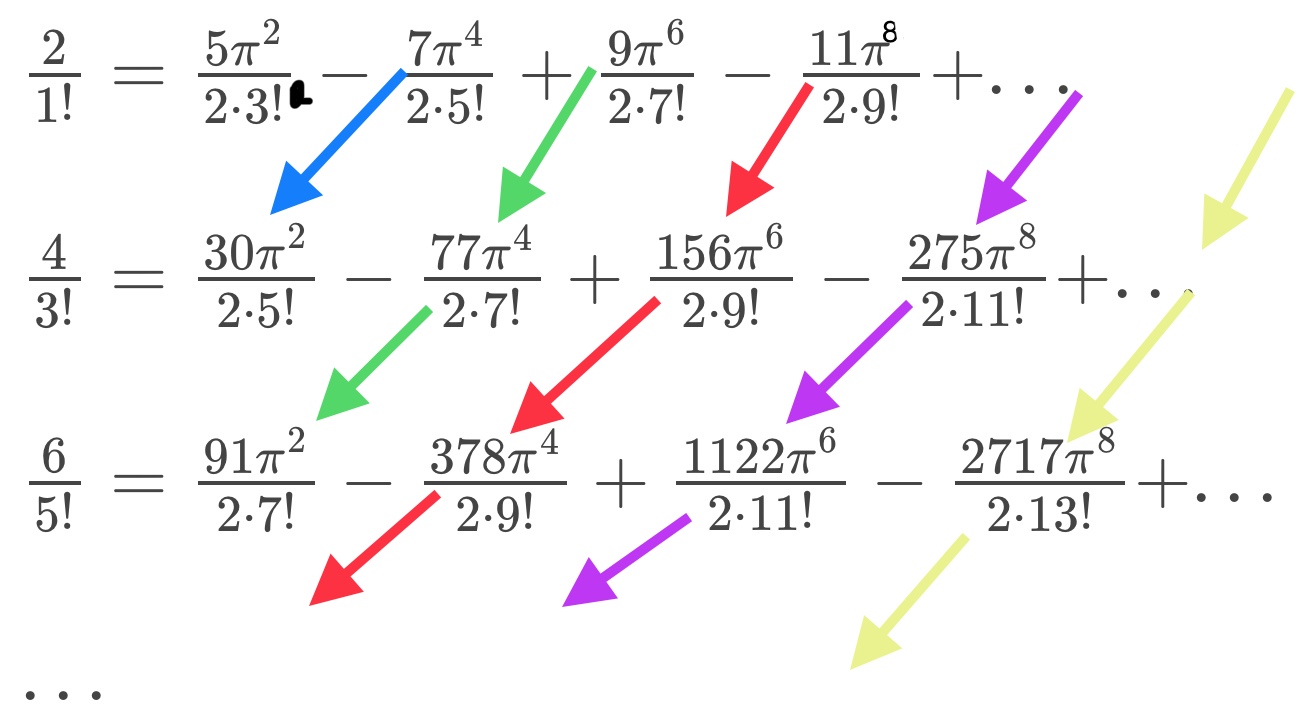}}
    \caption{Zig zag proof of Conjecture 6}
    \label{fig:my_label}
\end{figure}
Finally, here's a surprising relation between the other shallow diagonals of Pascal's triangle and its normal diagonals:
\\\medskip
\textbf{Conjecture 7.}
$$\frac{\pi}{2!}-\frac{\pi^3-2\pi}{4!}+\frac{\pi^5-4\pi^3+3\pi}{6!}-\frac{\pi^7-6\pi^5+10\pi^3-4\pi}{8!}+\cdots=2\Bigg(\frac{\frac{1}{0!}}{\pi}+\frac{\frac{2}{0!}}{\pi^3}+\frac{\frac{6}{0!}-\frac{\pi^2}{2!}}{\pi^5}+\frac{\frac{20}{0!}-\frac{4\pi^2}{2!}}{\pi^7}+\frac{\frac{70}{0!}-\frac{15\pi^2}{2!}+\frac{\pi^4}{4!}}{\pi^9}+\cdots\Bigg)$$
\\\medskip
\begin{figure}
    \centering
    \scalebox{0.14}
    {\includegraphics{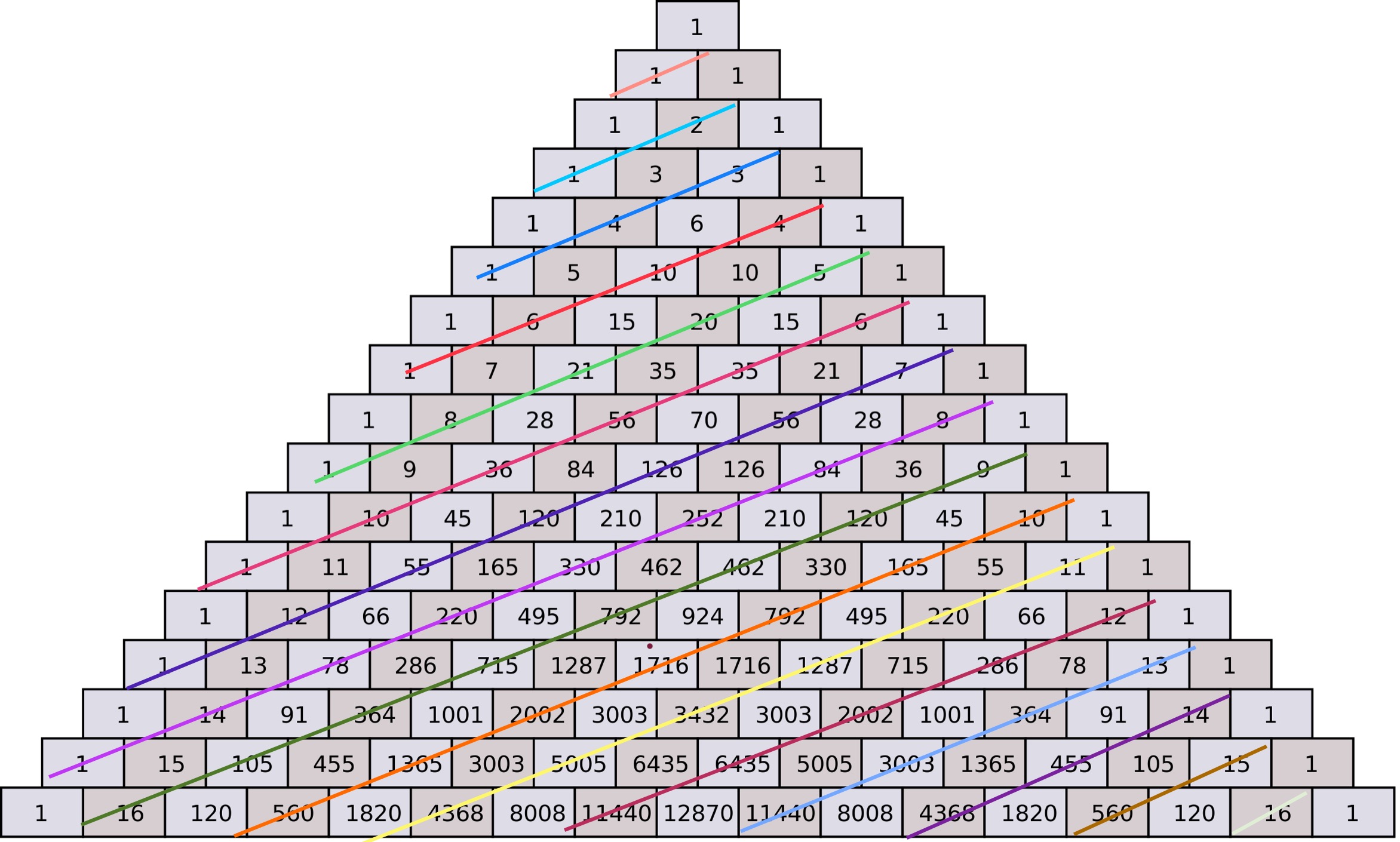}}
    \caption{LHS in Conjecture 7}
    \label{fig:my_label}
\end{figure}
\begin{figure}
    \centering
    \scalebox{0.3}
    {\includegraphics{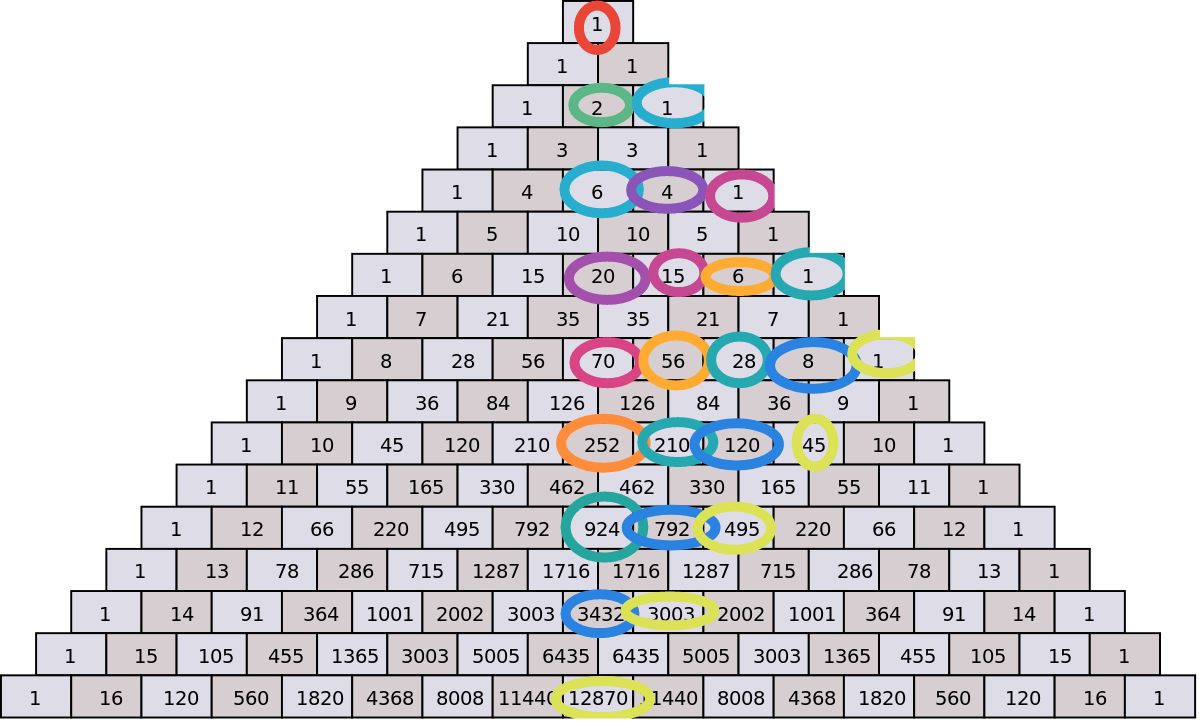}}
    \caption{RHS in Conjecture 7}
    \label{fig:my_label}
\end{figure}
\newpage
\textbf{REFERENCES}\\\medskip 
\medskip
[1] Harlan J. Brothers, Math bite: finding e in Pascal’s triangle, Math. Mag. 85 (1) (2012) 51.\\\medskip
[2] Gardner, Martin. “The Fibonacci sequence”. \textit{When You Were a Tadpole and I Was a Fish and Other Speculations About This and That}. First Edition, Hill and Wang, 2014, p. 113.\textit  {Archive.org}.  \url{https://archive.org/details/whenyouweretadpo00gard/page/n7/mode/2up}\\\medskip 
[3] Wolfram Alpha LLC. 2009. Wolfram|Alpha. \url{https://www.wolframalpha.com/input?i2d=true&i=Sum%5BDivide%5BDivide%5BPower%5B%CF%86%2Ck%2B1%5D-Power%5B%5C%2840%291-%CF%86%5C%2841%29%2Ck%2B1%5D%2Bx+Power%5B%CF%86%2Ck-1%5D-x+Power%5B%5C%2840%291-%CF%86%5C%2841%29%2Ck-1%5D+%2CSqrt%5B5%5D%5D%2Ck%21+%5D%2C%7Bk%2C0%2C%E2%88%9E%7D%5D&lang=es} 
(access December 16, 2022).\\\medskip 
[4] Weisstein, Eric W. "Fibonacci Polynomial." From MathWorld--A Wolfram Web Resource. \url{https://mathworld.wolfram.com/FibonacciPolynomial.html} 
\\\medskip
[5] Weisstein, Eric W. "Lucas Polynomial." From MathWorld--A Wolfram Web Resource. \url{https://mathworld.wolfram.com/LucasPolynomial.html}

\end{document}